\documentclass[10pt]{amsart}
\usepackage{amsmath}
\usepackage{cases}
\usepackage{mathrsfs}
\usepackage[utf8]{inputenc} 
\usepackage{CJKutf8}
\usepackage{fullpage}
\usepackage{bbm}
\usepackage{amssymb}
\usepackage{amscd}
\usepackage{indentfirst}
\usepackage{amsfonts,latexsym,amsthm,amsxtra,mathdots,latexsym,mathabx}
\usepackage[all,cmtip]{xy}
\RequirePackage{amsmath} \RequirePackage{amssymb}
\usepackage{color}
\usepackage{colordvi}
\usepackage[T1]{fontenc}
\usepackage{multicol}
\usepackage{hyperref}
\usepackage{mathtools}
\usepackage{comment}
\usepackage[margin=1.1in]{geometry}
\usepackage{xcolor}
\usepackage{parskip}
\allowdisplaybreaks[4]
\hypersetup{
	colorlinks,
	linkcolor={red!50!black},
	citecolor={blue!50!black},
	urlcolor={blue!80!black}
}

\newtheorem{thm}{Theorem}[section] 
\newtheorem{lem}[thm]{Lemma}

\newtheorem {rem}[thm]{Remark} 
\newtheorem* {notation}{Notation}

\newtheorem{lemma}[thm]{Lemma}

\newtheorem{theorem}[thm]{Theorem}

\newcommand{\Kl}{\operatorname{Kl}_3}
\newcommand{\Z}{\mathbb Z}

\newcommand{\e}{\mathrm e}
\newcommand{\eps}{\varepsilon}
\newcommand{\inv}[1]{\overline{#1}}
\newcommand{\sstar}{\mathop{\sum}\nolimits^{{\star}}}
\newcommand{\supp}{\operatorname{supp}}

\numberwithin{equation}{section}

\title{Mixed moment of twisted modular $L$-functions and Dirichlet $L$-functions}
\author{Zhaoyan Chen}
\date{}
\address{School of Mathematics, Shandong University, Jinan 250100, China}
\email{\href{mailto:zychen@mail.sdu.edu.cn}{zychen@mail.sdu.edu.cn}}

\subjclass[2020]{Primary 11M06; Secondary 11F12, 11L05.}
\keywords{Asymptotic formula, $L$-functions, character twists, Kloosterman sum}

\begin{document}
	\begin{abstract}
		Let $f$ be a Hecke holomorphic cusp form for $\mathrm{SL}(2,\mathbb{Z})$. Let $q=p^\gamma$ be a prime power, where $p$ is an arbitrary prime and $\gamma \geq 2$. We establish an asymptotic formula for the mixed moment of $L(s,f\otimes\chi)$ and $L(s,\chi)$ on the central point. We also establish an upper bound for the smoothed sum of Hecke eigenvalues twisted by hyper-Kloosterman sums modulo powers of $2$.
	\end{abstract}
	\maketitle
	
	\section{Introduction}

    Moments of central values of $L$-functions have long been studied due to its fundamental importance in number theory. A classical example is the $2k$-th moment of Dirichlet $L$-functions
	\[
	M_k(q):=\sideset{}{^{\star}}\sum_{\chi\bmod q}\left|L\!\left(\tfrac12,\chi\right)\right|^{2k},
	\]
	where $\sum^{\star}$ denotes the sum over primitive Dirichlet characters modulo $q$. Conrey et al.~\cite{CFKRS2005} conjectured that the general form for the $2k$-th moment is given by
	\[
	M_k(q) = \varphi^{{\star}}(q) P_{k^2}(\log q)+ O\!\left(q^{1/2+\varepsilon}\right),
	\]
	where $\varphi^{\star}(q)$ is the number of primitive characters modulo $q$, and $P_{k^2}(x)$ is a polynomial of degree $k^2$ whose coefficients depend on $q$. The above conjecture has been proved for $k=1$ (See \cite{IS1999} for example). For the case $k=2$, the asymptotic formula was first obtained by Heath-Brown~\cite{HB1981}. For $q\not\equiv 2\pmod 4$, he~\cite{HB1981} proved that
	\[
	\sideset{}{^{\star}}\sum_{\chi\bmod q}\left|L\!\left(\tfrac12,\chi\right)\right|^4= \frac{\varphi^{\star}(q)}{2\pi^2}\prod_{p\mid q}\frac{(1-p^{-1})^3}{1+p^{-1}}(\log q)^4+O\!\left(2^{\omega(q)}q(\log q)^3\right),
	\]
	where $\omega(q)$ denotes the number of distinct prime factors of $q$. However, its error term can exceed the main term when $q$ has many distinct prime factors. Subsequently, Soundararajan~\cite{Sound2007} overcame this difficulty and established the asymptotic formula with the same leading term for all $q\not\equiv 2\pmod 4$.
	
	However, obtaining a power saving error term remained a challenge even for prime moduli. It was only with the work of Young~\cite{Young2011} that such a saving was achieved. More precisely, he proved that
	\[
	\frac{1}{\varphi^{\star}(q)}\ 
	\sideset{}{^{\star}}\sum_{\chi\bmod q}
	\left|L\!\left(\tfrac12,\chi\right)\right|^4=P_4(\log q)
	+ O_{\varepsilon}\!\left(q^{-5/512+\varepsilon}\right)
	\]
	for prime $q$. Wu~\cite{Wu2023} subsequently extended Young's result to arbitrary moduli $q\not\equiv 2\pmod 4$. He obtained 
	\[
	\frac{1}{\varphi^\star(q)}\ \sideset{}{^{\star}}\sum_{\chi \bmod q}\left|L\!\left(\tfrac12,\chi\right)\right|^4=\prod_{p \mid q}\frac{(1-p^{-1})^3}{1+p^{-1}}P_4(\log q)+O\!\left(q^{-11/448+\varepsilon}\right).
	\]
	Despite many attempts, the asymptotic formulas $M_k(q)$ for higher moments have not been solved.
	
	The fourth moment of Dirichlet $L$-functions can also be viewed as a special case of a mixed moment of twists of Hecke eigenforms. For two fixed Hecke eigenforms $f$ and $g$, let us consider
	\[
	M_{f,g}(q)= \frac{1}{\varphi^{\star}(q)}\ \sideset{}{^{\star}}\sum_{\chi\bmod q}L\!\left(\tfrac12,f\otimes\chi\right)\overline{L\!\left(\tfrac12,g\otimes\chi\right)}.
	\]
	Let $E$ denote the Eisenstein series with normalized Hecke
	eigenvalues $d(n)$, where $d(n)$ is the divisor function.
	Since $L(s,E\otimes\chi)=L(s,\chi)^2$, taking $f=g=E$ recovers the fourth moment of Dirichlet $L$-functions. When both $f$ and $g$ are cusp forms, Blomer and Mili\'cevi\'c~\cite{BM2015} obtained an asymptotic formula with a power-saving error term for moduli $q$ satisfying a suitable factorization condition. More recently, Mili\'cevi\'c, Qin and Wu~\cite{MQW2025} and, independently, Pascadi~\cite{Pascadi2026} removed this restriction, obtaining power-saving asymptotics for arbitrary moduli $q\not\equiv 2\pmod 4$. Taking $f$ to be a cusp form and $g=E$ gives a cuspidal analogue of the fourth moment. For prime $q$, Blomer et al.~\cite{BFKMM2017} proved that
	\[
	M_{f,E}(q)
	= \frac{1}{\varphi^{\star}(q)}\ \sideset{}{^{\star}} \sum_{\chi\bmod q}
	L\!\left(\tfrac12,f\otimes\chi\right) \overline{L\!\left(\tfrac12,\chi\right)}^{\,2}= \frac{L(1,f)^2}{\zeta(2)}+ O_{f,\varepsilon}\!\left(q^{-1/68+\varepsilon}\right).
	\]
	The error term was subsequently improved to $O(p^{-1/64+\varepsilon})$ by Shparlinski~\cite{Shparlinski2019}, and further to $O(p^{-1/22+\varepsilon})$ for holomorphic $f$ and $O(p^{-5/152+\varepsilon})$ for Maass forms by Khan and Zhang~\cite{KZ2023}. Tang and Wu~\cite{TW2025} recently extended the result to general moduli, obtaining an error term comparable to the best bound in \cite{KZ2023} for prime moduli.
	
	Das and Khan~\cite{DK2015} studied a related moment for prime moduli. For a fixed Hecke--Maass cusp form $f$ for $\mathrm{SL}_2(\mathbb Z)$ and a prime $q$, they proved that
	\[
	\frac{1}{\varphi^{\star}(q)}\ \sideset{}{^{\star}}\sum_{\substack{\chi\bmod q\\ \chi(-1)=1}}L\!\left(\tfrac12,f\otimes\chi\right) \overline{L\!\left(\tfrac12,\chi\right)}=\frac{1}{2}\,L(1,f)
	+O_{g,\varepsilon}\!\left(q^{-8/64+\theta+\varepsilon}\right),
	\]
	where $\theta$ denotes the best bound towards the Ramanujan–Petersson conjecture for $f$, and one may take $\theta=7/64$ by Kim and Sarnak \cite{KS2003}. Subsequently, Sono~\cite{Sono2019} extended this mixed moment to points $s_0$ in the strip $\tfrac12\leq\Re(s_0)<1$.
	
	One may also ask how the central values of two different $L$-functions behave when twisted by the same character. Liu~\cite{Liu2015} studied the mixed moment of $L(\tfrac12,\pi\otimes\chi)L(\tfrac12,\chi)$, where $\pi$ is a fixed unitary cuspidal automorphic representation of $\mathrm{GL}_2(\mathbb A_{\mathbb Q})$. Set
	\[
	\mathcal D(Q)= \bigl\{ pr :Q^{3/4}<p\leq 2Q^{3/4},\ Q^{1/4}<r\leq 2Q^{1/4},\ p,r \text{ prime} \bigr\}.
	\]
	Liu~\cite{Liu2015} proved that
	\[
	\sum_{q\in\mathcal D(Q)}\ \sideset{}{^{\star}}\sum_{\substack{\chi\bmod q\\ \chi(-1)=1}}L\!\left(\tfrac12,\pi\otimes\chi\right)L\!\left(\tfrac12,\chi\right)=\frac12\sum_{q\in\mathcal D(Q)}q+ O_{\pi,\varepsilon}\!\left(Q^{15/8+\varepsilon}\right).
	\]

	For the two moments considered by Liu~\cite{Liu2015} and by Das and Khan~\cite{DK2015}, the effect of conjugation can be seen directly in the dual sum. Both moments detect simultaneous nonvanishing, but the existance of conjugation substantially simplifies the dual sum. Let $\epsilon(\chi)=\tau(\chi)/\sqrt q$ denote the normalized Gauss sum. Then the Gauss sum factor is $\epsilon(\chi)^3$ in Liu's moment~\cite{Liu2015} and $\epsilon(\chi)^2\overline{\epsilon(\chi)}=\epsilon(\chi)$ in the moment of Das and Khan~\cite{DK2015}. Thus, averaging over characters gives rise to a hyper-Kloosterman sum \(\mathrm{Kl}_3\) in the former case, while in the latter the corresponding $\mathrm{Kl}_1$ reduces to a single additive character.
	
	Motivated by these two results, we ask whether Liu's unconjugated moment~\cite{Liu2015} can be evaluated for a prime-power modulus $q$. Let $f$ be a fixed normalized holomorphic Hecke cusp form for $\mathrm{SL}_2(\mathbb Z)$, and let $q=p^\gamma$, where $p$ is a fixed prime. We consider
	\[
	M_f(q) = \sideset{}{^{\star}}\sum_{\substack{\chi\bmod q\\ \chi(-1)=1}} L\!\left(\tfrac12,f\otimes\chi\right)L\!\left(\tfrac12,\chi\right).
	\]
	The unconjugated moment is of interest in its own right. Since the same character twists both factors, $M_f(q)$ is the first moment of the Dirichlet twists of the fixed degree-three
	$L$-function $L(s,f)\zeta(s)$. In the conjugated moment, the two factors are twisted by $\chi$ and $\bar\chi$, respectively. Theorem~\ref{main1} thus determines the mean central value of this degree-three twist family for the prime-power modulus $p^\gamma$.
	\begin{theorem}\label{main1}
		Let $f$ be a Hecke holomorphic cusp form for $\mathrm{SL}(2,\mathbb{Z})$, and let $q=p^\gamma$ be a prime power with $\gamma\geq2$. Then we have
		\begin{equation}
			\sideset{}{^{\star}}\sum_{\substack{\chi \bmod q \\ \chi(-1)=1}}L\left(\frac{1}{2}, f \otimes \chi\right) {L\left(\frac{1}{2}, \chi\right)}=\frac{(p-1)^2 p^{\gamma-2}}{2}+O_p\left(q^{1-\varepsilon}\right).
		\end{equation}
	\end{theorem}
        In order to prove Theorem~\ref{main1}, we will use a power-saving estimate for an average of $\mathrm{GL}(2)$ Fourier coefficients twisted by hyper-Kloosterman sums. Define
        \[
        T(a,M;q):=\sum_{m\ge1}\lambda_f(m)\Kl(am;q)V\left(\frac{m}{M}\right),
        \]
        where $M$ is a large parameter and $V \in C_c^\infty((0,\infty))$ satisfies $V^{(k)} \ll q^{k\varepsilon}$ for all integers $k \geq 0$. The case of odd primes was proved by Sharma~\cite{Sharma2024}. For $p=2$, the derivative $2x$ of $x^2-t$ is even, so the simple-root version of Hensel's lemma does not apply. The square-root congruences arising in the proof therefore require a separate analysis, which allows us to extend the result to $p=2$ as follows:
		\begin{theorem}\label{main2}
			Let $\gamma\ge2$, $q=2^\gamma$, $(a,2)=1$. If $M\geq1$ and 
			\begin{equation}\label{range}
				1+M/q\le q^{1/10},
			\end{equation}
			then, for any $\varepsilon>0$,
			\begin{equation}\label{strong}
				T(a,M;q)\ll q^{1/3+\eps}M^{1/2}(1+M/q)^{2/3}.
			\end{equation}
			The implied constant is uniform in odd $a$, $M$, and $\gamma$.
		\end{theorem}
        \begin{notation}
        For $x,y,z\in\mathbb{Z}$, we denote by
		\[
			\Kl(a;q)=\frac1q\ \sideset{}{^{\star}}\sum_{x,y\bmod q}e_q(x+y+a\inv{xy}),
			\qquad S(1,b;q)=\sideset{}{^{\star}}\sum_{z\bmod q}e_q(z+b\bar z).
		\]
		As usual, $e(z)=e^{2\pi iz}$ and we write $e_Q(z)=e(z/Q)$. The ${\star}$ in $\sideset{}{^{\star}}\sum_{a\bmod q}$ means that the sum is restricted to $(a,q)=1$. Throughout the paper, $\varepsilon$ denotes an arbitrarily small positive constant, not necessarily the same at each occurrence.
        \end{notation}
        
		\section{Preliminaries}
		\subsection{Approximate functional equations} 
        Let $f$ be a holomorphic cusp form of weight $\kappa$ for $\mathrm{SL}(2,\mathbb{Z})$. Let $\lambda_{f}(n)$ denote the eigenvalue of the $n$-th Hecke operator corresponding to $f$. For an even Dirichlet character $\chi$ modulo $q$, we define the $L$-functions
        \[
        L(s,\chi)=\sum_{n=1}^{\infty}\chi(n)n^{-s} \qquad\text{and}\qquad L(s,f\otimes\chi)=\sum_{n=1}^{\infty}\lambda_f(n)\chi(n)n^{-s}.
        \]
		for $\operatorname{Re}(s)>1$, which admit analytic continuation to entire functions. We set
        \[
        \Lambda(s,\chi)=\left(\frac{q}{\pi}\right)^{s/2}\Gamma\left(\frac{s}{2}\right)L(s,\chi),\qquad\Lambda(s,f\otimes\chi):=\left(\frac{q}{2\pi}\right)^s\Gamma\left(s+\frac{\kappa-1}{2}\right)L(s,f\otimes\chi).
        \]
        and then the $L$-function satisfies the functional equations
        \[
        \Lambda(s,\chi)=\frac{\tau(\chi)}{\sqrt{q}}\Lambda(1-s,\overline{\chi}),\qquad\Lambda(s,f\otimes\chi)=\frac{\tau(\chi)^2}{q}\Lambda(1-s,f\otimes\overline{\chi}).
        \]
        The following lemma can be derived from \cite[Theorem 5.3]{IK2004}.
		\begin{lem}
			We have
			\begin{equation}\label{afe}
				L(\frac{1}{2}, f\otimes\chi)L(\frac{1}{2}, \chi)=\sum_{m,n}\frac{\lambda_f(n)\chi(mn)}{\sqrt{mn}}W\left(\frac{mn}{q\sqrt{q}}\right)+\epsilon_{\chi}^{3}\sum_{m,n}\frac{\lambda_f(n)\overline{\chi}(mn)}{\sqrt{mn}}W\left(\frac{mn}{q\sqrt{q}}\right),
			\end{equation}
			where for $y>0$,
            \begin{equation}\label{W}
			W(y)=\frac{1}{2 \pi i} \int_{(2)} \frac{\Gamma(\kappa/ 2+s)\Gamma(1/ 4+s/2)}{\Gamma(\kappa/ 2)\Gamma(1 / 4)} G(s) y^{-s} \frac{\mathrm{~d} s}{s}.
			\end{equation}
            $G(u)$ be any function which is holomorphic and bounded in the strip $-4<\operatorname{Re(u)}<4$, even and normalized by $G(0)=1$.
		\end{lem}
		\begin{rem}
			We have that $W(n)\ll\big(\frac{n}{q^{3/2}}\big)^{-A}$ for any $A>0$, thus the sums in \eqref{afe} are essentially supported on $mn<q^{3/2+\varepsilon}$. 
		\end{rem}
		
		\subsection{The character sums} 
		\begin{lem}\label{character}
			Let $p$ be a prime, let $\gamma\ge2$, and put $q=p^\gamma$. If $(n,q)=1$, then we have
			$$
			\begin{aligned}
				\sideset{}{^{\star}}\sum_{\substack{\chi \bmod q\\\chi(-1)=1} }\chi(n)=\frac{\varphi(q)}{2}(\delta_{n \equiv 1 \bmod q}+\delta_{n \equiv-1 \bmod q})-\frac{\varphi(q/p)}{2}(\delta_{n \equiv 1 \bmod q/p}+\delta_{n \equiv-1 \bmod q/p})
			\end{aligned}
			$$
			and
			$$
			\sideset{}{^{\star}}\sum_{\substack{\chi \bmod q\\\chi(-1)=1} } \epsilon(\chi)^3 \bar{\chi}(n)=\frac{\varphi(q)}{2 \sqrt{q}} (\mathrm{Kl}_3(n ; q)+\mathrm{Kl}_3(-n ; q)).
			$$
		\end{lem}
        \begin{proof}
        The first formula follows from the identity that 
        \[
            \sideset{}{^{\star}}\sum_{\substack{\chi \bmod q\\ \chi(-1)=1}} \chi(n)=\frac12\sum_{\pm}\sum_{\substack{q=dr\\ r\mid n\pm 1}}\mu(d)\varphi(r).
        \]
        Recall that, when $\gamma\geq 2$, $\tau(\chi)=0$ unless $\chi$ is primitive. We may therefore remove the restriction to primitive characters and apply the orthogonality of Dirichlet characters to obtain
        \[
       \begin{aligned}
        \sideset{}{^{\star}}\sum_{\substack{\chi \bmod q\\\chi(-1)=1} } \epsilon(\chi)^3 \bar{\chi}(n)&=\frac{1}{2q^{3/2}}\sum_{\chi \bmod q}(1+\chi(-1))\mathop{\sstar\sstar\sstar}_{x,y,z\bmod q}\e_q(x+y+z) \chi(xyz)\bar{\chi}(n)\\
        &=\frac{1}{2q^{3/2}}\mathop{\sstar\sstar\sstar}_{x,y,z\bmod q}\e_q(x+y+z) \left(\sum_{\chi \bmod q}\chi(xyz)\bar{\chi}(n)+\sum_{\chi \bmod q}\chi(-xyz)\bar{\chi}(n)\right)\\
        &=\frac{\varphi(q)}{2 \sqrt{q}} \left(\Kl(n ; q)+\Kl(-n ; q)\right).
        \end{aligned}
        \]
        This establishes the second formula.
        \end{proof}
		
		\subsection{Sums of Fourier coefficients}
        By the Rankin-Selberg theory, we have
	    \begin{lem} 
        We have 
		$$
		\begin{aligned}
			\sum_{n\leq x}|\lambda_{f}(n)|^{2}\ll x .
		\end{aligned}
		$$
	      \end{lem}
        For odd primes, we have the following result corresponding to Theorem \ref{main2}.
        \begin{lem}\label{oddprime}
        Let $q=p^{\gamma}$ with $p$ an odd prime and $\gamma\geq 2$. Assume $(a,q)=1$, we have
        \[
        \sum_{m\geq 1}\lambda_f(m)\mathrm{Kl}_3(am;q)V\left(\frac{m}{M}\right)\ll p^{7/12}M^{1/2}q^{1/3+\varepsilon}\left(1+\frac{M}{q}\right)^{2/3}+q^{13/20+\varepsilon}.
        \]
        \end{lem}
        \begin{proof}
        See \cite[Theorem 1.3]{Sharma2024}, taking $\mathcal{N}=\{1\}$ and $\alpha_1=1$.
        \end{proof}
        
		We will also need the Voronoi summation formula for the holomorphic cusp form. For $h(x)\in C_C(0,\infty)$, we set
		\[
		\Phi_h(x)=2\pi i^\kappa\int_0^\infty h(y)J_{\kappa-1}\left(4\pi\sqrt{xy}\right)\mathrm{d} y,
		\]
		where $J_{\kappa-1}$ is the usual $J$-Bessel function of order $\kappa-1$. We have the following Voronoi summation formula (see \cite[Theorem A.4]{KMV2002}).
		
		\begin{lemma}\label{Vor1}
			Let $q\in\mathbb{N}$ and $a\in\mathbb{Z}$ with $(a,q)=1$. For $X>0$, we have
			\[
			\sum_{n=1}^{\infty}\lambda_f(n)e\left(\frac{an}{q}\right)h\left(\frac{n}{X}\right)=\frac{X}{q}\sum_{n=1}^{\infty}\lambda_f(n)e\left(-\frac{\bar{a}n}{q}\right)\Phi_h\left(\frac{nX}{q^2}\right),
			\]
			where $\bar{a}$ denotes the multiplicative inverse of $a$ modulo $q$.
		\end{lemma}
		The function $\Phi_h(x)$ has the following properties (see \cite[Lemma 3.1]{LS2021}).
		\begin{lemma}
			For any fixed integer $J\ge 1$ and $x\gg 1$, we have
			\[
			\Phi_h(x)=x^{-1/4}\int_0^\infty h(y)y^{-1/4}\sum_{j=0}^{J}\frac{c_j e(2\sqrt{xy})+d_j e(-2\sqrt{xy})}{(xy)^{j/2}}\mathrm{d} y+O_{\kappa,J}\left(x^{-J/2-3/4}\right),
			\]
			where $c_j$ and $d_j$ are constants depending on $\kappa$.
		\end{lemma}
		For $x\ll1$, usually the bound $x^jJ_{\kappa-1}(x)\ll1$ would be sufficient in the applications.
		\subsection{Estimation of integrals}
		We will use the following stationary phase lemma several times.
		\begin{lem}\label{stationary} Suppose $w$ is a smooth function with compact support on $[Z, 2 Z]$ such that $w^{(j)}(\xi) \ll(Z / X)^{-j}$. Suppose that on the support of $w$, $h$ is smooth and satisfies that $h^{(j)}(\xi) \ll \frac{Y}{Z^{j}}$, for all $j \geq 0$. Let
			
			$$
			I=\int_{\mathbb{R}} w(\xi) e^{i h(\xi)} \mathrm{d} \xi.
			$$
			(i) If $h^{\prime}(\xi) \gg \frac{Y}{Z}$ for all $\xi \in \operatorname{supp} w$. Suppose $Y / X \geq 1$. Then $I \ll_A Z(Y / X)^{-A}$ for $A$ arbitrarily large.\\
			(ii) If $h^{\prime \prime}(\xi) \gg \frac{Y}{Z^2}$ for all $\xi \in \operatorname{supp} w$, and there exists $\xi_0 \in \mathbb{R}$ such that $h^{\prime}\left(\xi_0\right)=0$. Suppose that $Y / X^2 \geq R\geq1$. Then we have
			
			$$
			I=\frac{e^{i h\left(\xi_0\right)}}{\sqrt{h^{\prime \prime}\left(\xi_0\right)}} W(\xi_0)+O_A(ZR^{-A}),\quad \text { for any } A>0,
			$$
			for some smooth functions $W(\xi)$ (depending on $A$) supported on $\xi\asymp Z$.
		\end{lem}
        \begin{proof}
        See \cite[Lemma 3.1]{KPY2019}.
        \end{proof}

        \subsection{p-adic exponent pair}
		\begin{lemma}\label{klformula}
			Let $q=2^\gamma$ with $\gamma\ge2$. If $b$ is even, then $\Kl(b;q)=0$. If $b$ is odd, let $v$ be the unique odd cube root modulo $q$, so that $v^3\equiv b\pmod q$. Then
			\begin{equation}\label{klvalue}
				\Kl(b;q)=(-1)^\gamma e_q(3v).
			\end{equation}
			In particular, $|\Kl(b;q)|=1$ for odd $b$.
		\end{lemma}
		\begin{proof}
			For convenience, we write
			\begin{equation}
				\label{Kl}\Phi_b(x,y)=x+y+b\overline{xy},
			\end{equation}
			If $b$ is even, then shifting $x\mapsto x+q/2$ in \eqref{Kl} gives
			$$
			\begin{aligned}
				\Phi_b(x+q/2,y)&\equiv x+y+b\overline{xy}+\frac q2-\frac{bq\overline{x^2 y}}{2}\pmod q\\&\equiv x+y+b\overline{xy}+\frac q2\pmod q.
			\end{aligned}
			$$
			Thus each summand is multiplied by $e_q(q/2)=-1$, so the sum is equal to its negative and hence vanishes. 
			
			We now treat the cases where $b$ is odd and even separately. Since the even case is similar and simpler, we only give the details for odd $b$. Suppose $\gamma=2j+1$,we make a change of variable $x=x_0+2^{j+1}x_1,y=y_0+2^{j+1}y_1$. Note that 
            \begin{equation}\label{bar}
            \overline{u+p^k t}\equiv \overline{u}-\overline{u}^2p^kt+\overline{u}^3p^{2k}t^2\bmod p^{3k},
            \end{equation}
            Applying \eqref{bar}, we have
			\[
			\begin{aligned}
				\Kl(b;q)&=\frac 1q\sideset{}{^{\star}}\sum_{x_0\bmod 2^{j+1}}\sideset{}{^{\star}}\sum_{y_0\bmod 2^{j+1}}\sum_{x_1\bmod 2^j}\sum_{y_1\bmod 2^j}e_q(x_0+y_0+b\overline{x_0y_0})e_{2^j}(x_1(1-b\overline{{x_0}^2 y_0}))e_{2^j}(y_1(1-b\overline{x_0{y_0}^2}))\\&=\frac 12\sideset{}{^{\star}}\sum_{\substack{x_0\bmod 2^{j+1}}}
				\sideset{}{^{\star}}\sum_{\substack{y_0\bmod 2^{j+1}\\ y_0\equiv x_0\equiv v(\bmod 2^j)}}e_q(x_0+y_0+b\overline{x_0y_0}).
			\end{aligned}
			\]
			Again making the change of variables $x_0=v+2^j s$ and $y_0=v+2^jt$, and applying \eqref{bar}, we obtain
			\[
			\begin{aligned}
				\Kl(b;q)&=\frac 12\sum_{s\bmod 2}\sum_{t\bmod 2}e_q(3v)e_q(2^{2j}\bar{v}(s^2+st+t^2))\\&=\frac{e_q(3v)}{2}\sum_{s\bmod 2}\sum_{t\bmod 2}(-1)^{s^2+st+t^2}.
			\end{aligned}
			\]
			Since the sum over $s$ and $t$ is $-2$, we have $\Kl(b;q)=-e_q(3v)$.
		\end{proof}
        The following result can be derived using the method of $p$-adic exponent pairs (see \cite[Chapter 3]{Milicevic2016}).
		\begin{lem}\label{padic2}
			Let $q=2^\gamma$ and assume that $m$ is odd. For
			$V\in C_c^\infty((0,\infty))$, define
			\[
			S(m,L)=\sum_{l\ge1}\mathrm{Kl}_3(ml;q)V\left(\frac lL\right).
			\]
			Suppose that $(k,\ell,r,\delta,(n_0,u_0,\kappa_0,\lambda_0))$ is a $p$-adic exponent datum as in \cite[Definition 2]{Milicevic2016}. Then, for every $\kappa\ge\max\{2,\kappa_0(2/3,2)\}$ with $r=r(2/3,2,\kappa,\infty)$, assume that
            \[
            \gamma\geq\max\{2\kappa+2,\kappa+n_0(2/3,2,\kappa,\infty)\},\qquad 2^\kappa\le L\le q2^{-\kappa-2}.
            \]
            Then we have
			\[
			S(m,L)\ll2^{r+(1-k-\ell)\kappa-k}q^kL^{\ell-k}(\log q)^\delta.
			\]
		\end{lem}
		
		\begin{proof}
			We make a change of variable $l=c+2^\kappa t$, where $c$ is odd. Define 
			\[
			F_c(t)=3(mc)^{1/3}\left\{(1+2^\kappa\bar{c}t)^{1/3}-1\right\}.
			\]
			By uniqueness of the cube root, we have
			\[
			\bigl(m(c+2^\kappa t)\bigr)^{1/3}=(mc)^{1/3}(1+2^\kappa\bar{c}t)^{1/3}.
			\]
			Applying \eqref{klvalue}, we obtain
			\[
			S(m,L)=(-1)^\gamma\sum_{\substack{1\le c<2^\kappa\\c\ \mathrm{odd}}}e_q(3(mc)^{1/3})\sum_{t\in\mathbb Z}e_q(F_c(t))V\left(\frac{c+2^\kappa t}{L}\right).
			\]
			We next check that the above sum satisfies the hypotheses of the exponent-pair method. Since $1/3\in\mathbb Z_2$, we may expand
			\[
			(1+2^\kappa\bar{c}t)^{1/3}=\sum_{j\ge0}\binom{1/3}{j}2^{\kappa j}\overline{c^j}t^j
			\]
			whose coefficients lie in $\mathbb Z_2$ and tend to zero $2$-adically. Thus, $F_c$ is a convergent power series of the required type. Differentiating, we obtain
			\[ 
			F_c'(t)=2^\kappa\bar{c}(mc)^{1/3}(1+2^\kappa\bar{c}t)^{-2/3}.
			\]
			This is the defining form of the phase class with
			\[
			w=\kappa,\qquad y=\frac23,\qquad\omega=\bar{c},\qquad \omega'=\bar{c}(mc)^{1/3},
			\]
			and with the additional constant and error terms both equal to zero. In particular, the error-precision
			parameter $u$ may be chosen sufficiently large to satisfy the corresponding $u_0$ condition. Since
			\[
			\iota(2/3)=0,\qquad \iota'(2/3)=1,
			\]
			the effective length in the exponent-datum estimate is
			\[
			H_0=2^{\gamma-w-\kappa-\iota'(2/3)}=2^{\gamma-2\kappa-1}.
			\]
			The assumptions on $\kappa,\lambda,\gamma$ ensure that all remaining conditions of the datum are satisfied. Hence, for every integer $T$ and every $1\le B\le H_0$, applying Definition~2 of \cite{Milicevic2016}, we have
			\[
			\left|\sum_{T<t\le T+B}e_q(F_c(t))\right|\ll2^{r}\left(\frac{H_0}{B}\right)^kB^\ell(\log H_0)^\delta\ll2^{r}H_0^kB^{\ell-k}(\log q)^\delta.
			\]
			Put $H=L/2^\kappa$. The support of $V$ is contained in the interval
			\[
			\alpha_c<t<\beta_c,
			\qquad
			\alpha_c=\frac{L-c}{2^\kappa},
			\quad
			\beta_c=\frac{2L-c}{2^\kappa},
			\]
			whose length is $H$. Taking $T_c=[\alpha_c]$, every partial sum starting at $T_c$ and ending in this interval has integer length at most $H+1\le2H$. Moreover,
			\[
			2H=\frac{2L}{2^\kappa}\le 2^{\gamma-2\kappa-1}=H_0.
			\]
			Therefore
			\[
			\sup_{\alpha_c\le x\le\beta_c}\left|\sum_{T_c\le t\le x}e_q(F_c(t))\right|\ll2^{r}H_0^kH^{\ell-k}(\log q)^\delta.
			\]
			By partial summation, it follows that
			$$
			\begin{aligned}
				\sum_{t\in\mathbb Z}e_q(F_c(t))V\left(\frac{c+2^\kappa t}{L}\right)&=-\frac{2^\kappa}{L}\int_{\mathbb{R}}\left(\sum_{T_c\le t\le x}e_q(F_c(t))\right) V'\left(\frac{c+2^\kappa x}{L}\right)\mathrm{~d} x\\&\ll2^{r}H_0^kH^{\ell-k}(\log q)^\delta.
			\end{aligned}
			$$
			Summing over the odd residue classes modulo $2^\kappa$, we conclude that
			\[
			\begin{aligned}
				|S(m,L)|&\ll2^\kappa 2^{r}\left(2^{\gamma-2\kappa-1}\right)^k\left(\frac L{2^\kappa}\right)^{\ell-k}(\log q)^\delta\\&=2^{r+(1-k-\ell)\kappa-k}q^kL^{\ell-k}(\log q)^\delta.
			\end{aligned}
			\]
			This completes the proof.
		\end{proof}
        For the case of odd primes, see \cite[Lemma 3.2]{WZ2026}, where the following result is established.
        \begin{lem}\label{padicodd}
       With the notation as above, let $\kappa \geq \max\{\kappa_0(2/3,p),1+\iota'(2)\}$ and put $r=r(2/3,p,\kappa,\infty)$. Assume that
        \[
        \gamma \geq \max\{n_0(2/3,p,\kappa,\infty)+\kappa,2\kappa\},\qquad p^\kappa \leq L \leq qp^{-\kappa},
        \]
        when $p \geq 5$. When $p=3$, assume that
        \[
        \gamma \geq\max\{5,n_0(2/3,3,\kappa,\infty)+\kappa+1,2\kappa+1\},\qquad 3^{\kappa+1} \leq L \leq q3^{-\kappa}.
        \]
        Then
        \[
        S(m,L) \ll p^{r+(1-k-\ell)\kappa}q^k L^{\ell-k}(\log q)^\delta.
        \]
        \end{lem}
        
		\section{Proof of Theorem \ref{main1}}
		Applying the approximate functional equation \eqref{afe} and Lemma \ref{character}, we may write 
		\[
		\sideset{}{^{\star}}\sum_{\substack{\chi \bmod q\\\chi(-1)=1}} L\left(\frac{1}{2}, f \otimes \chi\right) L\left(\frac{1}{2}, \chi\right)=S_1+S_2,
		\]
		where
		$$
		\begin{aligned}
			&S_1=\frac{\varphi(q)}{2}\sum_{\substack{m, n\geq 1\\mn\equiv\pm1(\bmod q)}}\frac{\lambda_f(n)}{\sqrt{mn}}W\left(\frac{mn}{q\sqrt{q}}\right)-\frac{\varphi(q/p)}{2}\sum_{\substack{m, n\geq 1\\mn\equiv\pm1(\bmod q/p)}}\frac{\lambda_f(n)}{\sqrt{mn}}W\left(\frac{mn}{q\sqrt{q}}\right),\\
			&S_2=\frac{\varphi(q)}{2\sqrt{q}}\sum_{m,n\geq1}\frac{\lambda_f(n)}{\sqrt{mn}}\left(\mathrm{Kl}_3(mn; q)+\mathrm{Kl}_3(-mn; q)\right)W\left(\frac{mn}{q\sqrt{q}}\right).
		\end{aligned}
		$$
		We first consider $S_1$, for which we obtain the following lemma.
		\begin{lem}
			Given any $\varepsilon>0$, we have
			\[
			S_1=\frac{(p-1)^2 p^{\gamma-2}}{2}+O(q^{3/4+\varepsilon}).
			\]
		\end{lem}
		\begin{proof}
			The main contribution to $S_1$ arises from the terms with $mn=1$. By the definition of $V$, this contribution is given by
			$$
			\begin{aligned}
				\frac{1}{2 \pi i}\frac{\varphi(q)-\varphi(q/p)}{2} \int_{(2)} \frac{\Gamma(\kappa/ 2+s)\Gamma(1/ 4+s/2)}{\Gamma(\kappa/ 2)\Gamma(1 / 4)} G(s)q^{3s/2} \frac{\mathrm{~d} s}{s}
			\end{aligned}
			$$
			We now shift the contours of integration to $\operatorname{Re}(s)=-1/2+\varepsilon$, we pick up the residue $\frac{(p-1)^2 p^{\gamma-2}}{2}$ at $s=0$. The remaining integrals along the shifted contours can be bounded in a standard way by $O(q^{1/4+\varepsilon})$.
			
			It remains to consider the terms with $mn\neq 1$. For $r=q$ or $r=q/p$, the congruence
			\[
			mn\equiv \pm1\pmod r
			\]
			implies that
			\[
			mn=\ell r\pm1
			\]
			for some integer $\ell\geq1$. By the rapid decay of $V$, we may restrict to
			\[
			mn\ll q^{3/2+\varepsilon},
			\]
			and hence
			\[
			\ell\ll \frac{q^{3/2+\varepsilon}}{r}.
			\]
			Using Deligne's bound $\lambda_f(n)\ll_\varepsilon n^\varepsilon$ together with the divisor bound, we obtain
			$$
			\begin{aligned}
				\varphi(r)\sum_{\substack{mn\equiv\pm1(\bmod r)\\ mn\neq1}}\frac{|\lambda_f(n)|}{\sqrt{mn}}W\left(\frac{mn}{q^{3/2}}\right)&\ll rq^\varepsilon\sum_{\ell\ll q^{1/2+\varepsilon}}\frac{1}{\sqrt{\ell r}} \ll q^{3/4+\varepsilon}.
			\end{aligned}
			$$
			Combining the cases $r=q$ and $r=q/p$, we obtain the upper bound $O(q^{3/4+\varepsilon})$ for the terms in $S_1$ with $mn\neq 1$. This completes the proof.
		\end{proof}
		
		We now turn to $S_2$, for which we obtain the following lemma. 
		\begin{lem}
			Given any $\varepsilon>0$, we have
			\[
			S_2\ll q^{1-\varepsilon}.
			\]
		\end{lem}
		\begin{proof}
			It suffices to treat the term involving $\mathrm{Kl}_3(mn;q)$, since the term with $\mathrm{Kl}_3(-mn;q)$ can be handled in exactly the same way. We now apply a smooth partition of unity to the variables $m$ and $n$. Thus, $S_2$ is reduced to $O(\log^2 q)$ sums of the form
			\[
			\frac{\varphi(q)}{\sqrt{q}}\sum_{m,n\geq1}\mathrm{Kl}_3(mn; q)\frac{\lambda_f(n)}{\sqrt{mn}}W\left(\frac{mn}{q^{3/2}}\right)V_1\left(\frac{m}{M}\right)V_1\left(\frac{n}{N}\right),
			\]
			where $MN\ll q^{3/2+\varepsilon}$ and $V_1$ is a smooth function supported on $[1/2,5/2]$ satisfying
            \[
            \sideset{}{^d}\sum_M V_1\left(\frac{m}{M}\right)=1.
            \]
			Here $\sum_{M}^{d}$ denotes a sum over powers of $2$.
            
			\textbf{Case 1:} $MN\ll q^{1-\delta}$. In this range, a trivial estimate gives $S_2\ll \sqrt{qMN}\ll q^{1-\varepsilon}$ upon taking $\delta=2\varepsilon$.\\
			\textbf{Case 2:} $MN\gg q^{1+\delta}$. Denote by
			\[
			B(M,N)
			=
			\sum_{m,n\ge1}\mathrm{Kl}_3(mn;q)\frac{\lambda_f(n)}{\sqrt{mn}}W\left(\frac{mn}{q^{3/2}}\right)V_1\left(\frac mM\right)V_1\left(\frac nN\right).
			\]
			Since $\mathrm{Kl}_3(mn;q)=0$ when $p\mid mn$, we may restrict $m$ to $(m,q)=1$. Expanding the Kloosterman sum gives
			\[
			B(M,N)=\frac1q\sum_{\substack{m\ge1\\(m,q)=1}}\frac{1}{\sqrt{m}}\sideset{}{^{\star}}\sum_{x,y\bmod q}e_q(x+y)\sum_{n\ge1}\frac{\lambda_f(n)}{\sqrt{n}}e_q(m\overline{xy}n)W\left(\frac{mn}{q^{3/2}}\right)V_1\left(\frac mM\right)V_1\left(\frac nN\right).
			\]
            Recalling from \eqref{W} that
            \[
            W\left(\frac{mn}{q^{3/2}}\right)=\frac{1}{2 \pi i} \int_{(\varepsilon)} \frac{\Gamma(\kappa/ 2+s)\Gamma(1/ 4+s/2)}{\Gamma(\kappa/ 2)\Gamma(1 / 4)} \left(\frac{mn}{q^{3/2}}\right)^{-s} G(s) \frac{\mathrm{~d} s}{s},
            \]
			and the rapid decay of $\Gamma(s)$ (which follows from Stirling's formula), we see that it suffices to bound
			\[
			\frac1q\sum_{\substack{m\ge1\\(m,q)=1}}\frac{1}{m^{1/2+s}}\sideset{}{^{\star}}\sum_{x,y\bmod q}e_q(x+y)\sum_{n\ge1}\frac{\lambda_f(n)}{n^{1/2+s}}e_q(m\overline{xy}n)V_1\left(\frac mM\right)V_1\left(\frac nN\right)
			\]
			for $|s|\ll q^{\varepsilon}$ and $\operatorname{Re}(s)=\varepsilon$. We may absorb all the $s$-dependent factors, together with $\sqrt{\frac{MN}{mn}}$, into the smooth function $V_1$ without affecting its essential properties. Consequently, it suffices to estimate
			\begin{equation}\label{B}
				\frac{1}{q\sqrt{MN}}\sum_{\substack{m\ge1\\(m,q)=1}}\sideset{}{^{\star}}\sum_{x,y\bmod q}e_q(x+y)\sum_{n\ge1}\lambda_f(n)e_q(m\overline{xy}n)V_1\left(\frac mM\right)V_1\left(\frac nN\right),
			\end{equation}
			We first apply the Voronoi summation formula to the $n$-sum. This gives
			\[
			\sum_{n\ge1}\lambda_f(n)e_q(m\overline{xy}n)V_1\left(\frac nN\right)=\frac{N}{q}\sum_{n\ge1}\lambda_f(n)e_q(-\bar mxy n)\Phi_{V_1}\left(\frac{nN}{q^2}\right),
			\]
			where
			\[
			\Phi_{V_1}(x)=x^{-1 / 4} \int_0^{\infty} V_1(y) y^{-1 / 4} \sum_{j=0}^J \frac{c_j e(2 \sqrt{x y})+d_j e(-2 \sqrt{x y})}{(x y)^{j / 2}} \mathrm{~d} y+O\left(x^{-J / 2-3 / 4}\right).
			\]
			Note that $4\sqrt{xy}^{(i)}(y)\asymp \sqrt{x}$ for $y\asymp 1$. Hence, by Lemma \ref{stationary} (i), if $nN\gg q^{2+\varepsilon}$, then $\Phi_{V_1}\ll q^{-2026}$. Thus, we may restrict the $n$-sum to the range $nN\ll q^{2+\varepsilon}$. Substituting the Voronoi formula back into \eqref{B}, we arrive at
			\begin{equation}\label{T1}
				\frac{\sqrt{N}}{q^2\sqrt{M}}\sum_{\substack{m\ge1\\(m,q)=1}}V_1\left(\frac mM\right)\sum_{n\ge1}\lambda_f(n)\Phi_{V_1}\left(\frac{nN}{q^2}\right)\mathcal T(m,n),
			\end{equation}
			where
			\[
			\mathcal T(m,n)=\sideset{}{^{\star}}\sum_{x,y\bmod q}e_q(x+y-\bar mxy n).
			\]
			Summing first over $y$ gives
			\[
			\sideset{}{^{\star}}\sum_{x\bmod q}e_q(x)\left(q\,\mathbf1_{q\mid (1-\bar m xn)}-\frac qp\,\mathbf1_{q/p\mid (1-\bar m xn)}\right).
			\]
			If $p\mid n$, then it is easy to see that $\mathcal T(m,n)=0$. If $(n,q)=1$, then
			\[
			\begin{aligned}
				\mathcal T(m,n)&=q\,e_q(m\bar n)-\frac qp\sum_{\substack{x\bmod q\\x\equiv m\bar n\;(\bmod\,q/p)}}e_q(x)\\&=q\,e_q(m\bar n)-\frac{qe_q(m \bar{n})}{p}\sum_{j=0}^{p-1} e_p(j)\\&=q\,e_q(m\bar n),
			\end{aligned}
			\]
			Thus, \eqref{T1} becomes
			\[
			\frac{\sqrt{N}}{q\sqrt{M}}\sum_{\substack{n\ge1\\(n,q)=1}}\lambda_f(n)\Phi_{V_1}\left(\frac{nN}{q^2}\right)\sum_{\substack{m\ge1\\(m,q)=1}}e_q(m\bar n)V_1\left(\frac mM\right).
			\]
			
			Writing $m=a+hq$, with $a\bmod q$ satisfying $(a,q)=1$ and $h\in\mathbb Z$, we decompose the $m$-sum into sums over $a$ and $h$. Applying the Poisson summation formula to the $h$-sum, we obtain
			\[
			\begin{aligned}
				\sum_{\substack{m\ge1\\(m,q)=1}}e_q(m\bar n)V_1\left(\frac mM\right)&=\frac Mq\sum_{h\in\mathbb Z}\widehat V_1\left(\frac{hM}{q}\right)\sideset{}{^{\star}}\sum_{a\bmod q}e_q((\bar n+h)a)\\&=\frac Mq\sum_{h\in\mathbb Z}(q\,\mathbf1_{nh\equiv-1\;(\bmod\,q)}-\frac qp\,\mathbf1_{nh\equiv-1\;(\bmod\,q/p)})\widehat V_1\left(\frac{hM}{q}\right).
			\end{aligned}
			\]
			Consequently, \eqref{T1} can be rewritten as
			\begin{equation}\label{T2}
				\frac{\sqrt{MN}}{q^2}\sum_{\substack{n\ge1\\(n,q)=1}}\lambda_f(n)\Phi_{V_1}\left(\frac{nN}{q^2}\right) \sum_{h\in\mathbb Z}(q\,\mathbf1_{nh\equiv-1\;(\bmod\,q)}
				-\frac qp\,\mathbf1_{nh\equiv-1\;(\bmod\,q/p)}))\widehat V_1\left(\frac{hM}{q}\right).
			\end{equation}
			Since $\gamma\ge 2$, both $q$ and $q/p$ are greater than $1$. When $h=0$, neither of the congruences
			\[
			nh\equiv -1\pmod q,
			\qquad
			nh\equiv -1\pmod{q/p}
			\]
			can hold. Hence, the term corresponding to $h=0$ vanishes, and we may restrict the inner sum to $h\ne 0$. By Lemma \ref{stationary} (i), we may further restrict the $h$-sum to $h\ll q^{1+\varepsilon}/M$, with a negligible error. Recall that $n\ll q^{2+\varepsilon}/N$.

			For $d\in\{q,q/p\}$, the congruence $nh\equiv-1\pmod d$ is equivalent to $nh=jd-1$. There are $O\big(\frac{q^3}{NMd}\big)$ possible values of $j$, and $jd-1\ne0$ since $d\ge2$. For each $j$, the number of pairs $(n,h)$ is bounded by the divisor function $d(|jd-1|)$. This gives
			\[
			\sum_{\substack{n\ll q^{2+\varepsilon}/N,\;0<|h|\ll q^{1+\varepsilon}/M \\(n,q)=1,\;nh\equiv-1\;(\bmod\,d)}}|\lambda_f(n)|\ll_\rho q^\varepsilon\left(1+\frac{q^3}{NMd}\right).
			\]
			Therefore, we conclude that
			\[
			\begin{aligned}
				\eqref{T2}&\ll q^\varepsilon\frac{\sqrt{MN}}{q^2}\left[q\left(1+\frac{q^2}{MN}\right)+\frac qp\left(1+\frac{q^2 p}{MN}\right)\right]\\&\ll q^\varepsilon\frac{\sqrt{MN}}{q^2}\left(q+\frac{q^3}{MN}\right)\\&=q^\varepsilon\left(\frac{\sqrt{MN}}q+\frac{q}{\sqrt{MN}}\right).
			\end{aligned}
			\]
			
			Since $MN\le q^2$, it follows that
			\[
			\frac{\varphi(q)}{\sqrt q}|B(M,N)|\ll \frac{q^{3/2+\varepsilon}}{\sqrt{MN}}.
			\]
			In the present range $MN\gg q^{1+\delta}$, taking $\delta=4\varepsilon$, this is $O(q^{1-\varepsilon})$. After summing over the $O(\log^2 q)$ dyadic blocks, the contribution of this range is $O(q^{1-\varepsilon})$.
			
			\textbf{Case 3:} $q^{1-\delta_0}\ll MN\ll q^{1+\delta_0}$. As in Case 2, applying Stirling's formula to remove the weight function $W$, we are reduced to estimating
			\begin{equation}\label{case3}
				\frac{\varphi(q)}{\sqrt{qMN}}\sum_{m,n\geq1}\mathrm{Kl}_3(mn; q)\lambda_f(n)V_1\left(\frac{m}{M}\right)V_1\left(\frac{n}{N}\right).
			\end{equation}
			We treat separately the cases $p=2$ and $p\neq2$, beginning with the case $p=2$. In view of Theorem \ref{main2}, we have that
			\[
			\eqref{case3}\ll \sqrt{\frac{q}{MN}}q^{1/3+\eps}MN^{1/2}(1+N/q)^{2/3}\ll \sqrt{M}q^{5/6+2\delta_0/3+\varepsilon}.
			\]
			In the range $M\ll q^{1/3-2\delta_0}$, taking $\delta_0=6\varepsilon$, $\eqref{case3}$ is $O(q^{1-\varepsilon})$. For $q^{1/3-2\delta_0}\ll M\ll q$, we apply the method of $p$-adic exponent pairs. Using the $A$- and $B$-processes (see, for example, \cite[Theorem 2]{Milicevic2016}), we obtain the exponent pair $(k,\ell)=\left(\frac{11}{82},\frac{57}{82}\right)$. Lemma~\ref{padic2} then yields
			\[
			\eqref{case3}\ll q^{\delta_0/2}Nq^{11/82}M^{46/82}(\log q)^\delta\ll q^{81/82+\varepsilon}.
			\]
			For the remaining range of $M$, we apply the Poisson summation formula to the $m$-sum in \eqref{case3} and obtain
			\[
			\frac{M}{q}\sum_{m\in \mathbb{Z}}\sum_{a\bmod q}\Kl(an, q)e_{q}(am)\widehat V_1\left(\frac{mM}{q}\right).
			\]
			Expanding the Kloosterman sum gives
			\[
			\frac{M}{q}\sum_{m\in \mathbb{Z}}S(1,-\bar{m}n; q)\widehat V_1\left(\frac{mM}{q}\right).
			\]
			By the rapid decay of $\widehat V$, we may restrict the $m$-sum to $m\ll q^{1+\varepsilon}/M$. By the Weil bound for Kloosterman sums, \eqref{case3} is bounded by $O(q^{1/2+\varepsilon})$.

            For $p \neq 2$, Lemma \ref{oddprime} yields
            $$
            \begin{aligned}
            \eqref{case3}\ll \sqrt{\frac{q}{MN}}p^{\frac{7}{12}}q^{\frac{1}{3}+\varepsilon}M\sqrt{N}(1+N/q)^{\frac{2}{3}}+Mq^{\frac{13}{20}+\frac{1}{2}\delta_0+\varepsilon}\ll \sqrt{M}q^{\frac{5}{6}+\frac{7}{12\gamma}+\frac{2}{3}\delta_0+\varepsilon}+Mq^{\frac{13}{20}+\frac{1}{2}\delta_0+\varepsilon}.
            \end{aligned}
            $$
            In the range $M\ll q^{1/3-2\delta_0}$, choosing $\delta_0=12\varepsilon$ and $\gamma\geq \frac{7}{2\delta_0}$, \eqref{case3} can be bounded by $O(q^{1-\varepsilon})$. For $q^{1/3-2\delta_0}\ll M\ll q$, we further apply Lemma \ref{padicodd} with $(k,\ell)=\left(\frac{11}{82}, \frac{57}{82}\right)$ to obtain
            \[
            \eqref{case3}\ll q^{\frac{1}{2}\delta_0}p^{r+\frac{7}{41}\kappa}Nq^{\frac{11}{82}}M^{\frac{46}{82}}(\log q)^{\delta}\ll q^{\frac{81}{82}+\varepsilon},
            \]
            provided that $\gamma\geq\varepsilon^{-1}\left(r+\frac{7}{41}\kappa\right)$. And we can take $\kappa=1$ and $r=0$ when $p\notin\{3,5,7,11,13,17,19\}$. If $p$ belongs to the finitely many primes, then the above bound also holds for $\gamma\geq N_0$ with $N_0$ a fixed value, and the case $\gamma\leq N_0$ can be treated directly. For the remaining range of $M$, the same Poisson summation argument gives $O(q^{1/2+\varepsilon})$.
		\end{proof}
		
		\section{Proof of Theorem \ref{main2}}
		
		\subsection{Some Lemmas}
		\begin{lemma}\label{classical}
			Let $\gamma\ge8$, $q=2^\gamma$, and let $B\in\mathbb Z_2^\times$. If $B\not\equiv1\pmod8$, then $S(1,B;q)=0$. Otherwise, choose $w\in\Z_2^\times$ with $w^2=B$. Then
			\begin{equation}\label{classvalue}
				S(1,B;q)=c_\gamma(w)e_q(2w)+\overline{c_\gamma(w)}e_q(-2w),
			\end{equation}
			where
			\[
			c_{2j}(w)=2^j(1+i^{w^{-1}}),\qquad
			c_{2j+1}(w)=2^{j+1}e_8(w^{-1}).
			\]
			The inverses in the exponents need only be taken modulo $4$ and $8$, respectively. Hence $|c_\gamma(w)|=\sqrt{2q}$ and $|S(1,B;q)|\le2\sqrt{2q}$.
		\end{lemma}
		\begin{proof}
			We shall use use the stationary phase method to evaluate the Kloosterman sum. We make the change of variables $z=z_0+2^{\lceil\gamma/2\rceil}z_1$, which gives
			\[
			S(1,B;q)=\sideset{}{^{\star}}\sum_{z_0\bmod 2^{\lceil\gamma/2\rceil}}e_q(z_0+B\bar{z}_0)\sum_{z_1\bmod 2^{\lfloor\gamma/2\rfloor}}e_{2^{\lfloor\gamma/2\rfloor}}(z_1(1-B\overline{{z_0}^2})).
			\]
			The $z_1$-sum vanishes unless $1-B\overline{{z_0}^2}\equiv 0\pmod {2^{\lfloor\gamma/2\rfloor}}$. Since $z\equiv z_0\pmod{2^{\lfloor\gamma/2\rfloor}}$, we have $z^2\equiv B\pmod{2^{\lfloor\gamma/2\rfloor}}$. This proves the vanishing condition modulo $8$.
			
			If $\gamma=2j$ with $j\ge4$, we first write $z=z_0+2^{j}z_1$. It follows from the above argument that $z_0^2\equiv w^2\pmod {2^j}$. The four solutions modulo $2^j$ split into the
			two groups
			\[
			z_0\equiv w\pmod{2^{j-1}}
			\qquad\text{and}\qquad
			z_0\equiv-w\pmod{2^{j-1}}.
			\]
			The  first group are represented by $z_0=w+2^{j-1}l$ with $l\bmod 2$. The corresponding sum is
			\[
			2^je_q(2w)\sum_{l\bmod 2}e_4(\bar{w}l^2)=2^j(1+i^{\bar{w}})e_q(2w).
			\]
			The second group contribute the complex conjugate.
			
			If $\gamma=2j+1\ge9$, with $j\ge4$, instead take $l\bmod4$. Since $3(j-1)\ge2j+1$, the positive-root contribution is
			\[
			2^j e_q(2w)\sum_{l\bmod4}e_8(\bar{w}l^2)
			=2^{j+1}e_8(\bar{w})e_q(2w).
			\]
			Again, the negative-root contribution is its complex conjugate.
		\end{proof}
		
		Let $u\ge1$ and $h=\gamma-u\ge4$. Throughout this section, $s_j,t_j,b_j$ are odd. Define
		\begin{equation}\label{corrdef}
			\mathcal C(d)=\mathop{\sstar\sstar}_{\substack{x_1,x_2\bmod2^u\\b_1\bar x_1-b_2\bar x_2\equiv d(\bmod 2^u)}}
			S(1,(t_1+2^hs_1x_1)^{-1};q)
			\overline{S(1,(t_2+2^hs_2x_2)^{-1};q)}.
		\end{equation}
		The sum in \eqref{corrdef} has $2^{u-1}$ terms, each of absolute value at most $8q$. Hence, the trivial bound for $\mathcal C(d)$ is $O(2^{\gamma+u})$. Moreover, we have the following lemma, which gives a sharper estimate.
		\begin{lemma}\label{local}
			Assume $\gamma\ge8$ and $d\ne0$, and put $v=v_2(d)$. If $d$ is odd or $t_j\not\equiv1\pmod8$ for $j\in\{1,2\}$, then $\mathcal C(d)=0$. In the remaining cases, choose $r_j\in\mathbb Z_2^\times$ satisfying $r_j^2=t_j^{-1}$ in $\mathbb Z_2$.\\
			(1) If $h\ge\lfloor u/2\rfloor+1$ or $v\le h-3$, then
			\begin{equation}\label{corrbound}
				|\mathcal C(d)|\le32\,2^{\gamma+\lceil u/2\rceil+v}.
			\end{equation}
			(2) If $h\le\lfloor u/2\rfloor$ and $v\ge h-2$, then $\mathcal C(d)=0$ unless there exist signs $\epsilon_1,\epsilon_2\in\{\pm1\}$ such that
			\begin{equation}\label{corrsupport}
				\epsilon_1s_1r_1^3b_1\equiv\epsilon_2s_2r_2^3b_2\pmod{2^{h-2}}.
			\end{equation}
		\end{lemma}
		\begin{proof}
			For odd $d$, the constraint has no solution. The vanishing condition on $t_j$ follows from Lemma~\ref{classical}. Henceforth assume that $d$ is even and both $t_j\equiv1\pmod8$. The congruence condition on $x_1$ and $x_2$ gives
			\[
			x_2\equiv Y(x_1)=\frac{b_2x_1}{b_1-dx_1}\pmod{2^u}.
			\]
			Define $G_j(x)$ as follows, and we have that
			\[
			\begin{aligned}
				G_j(x)&:=2^{1-h}r_j\{(1+2^hs_jt_j^{-1}x)^{-1/2}-1\}\\&=2^{1-h}r_j\sum_{i\geq1}\binom{-\frac{1}{2}}{i}(2^hs_jt_j^{-1}x)^i=2^{1-h}r_j\sum_{i\geq1}(-1)^i2^{-2i}\binom{2i}{i}(2^hs_jt_j^{-1}x)^i
			\end{aligned}
			\]
			The coefficient of $x^i$ has valuation $(h-2)i+1-h+s_2(i)$, which is nonnegative and tends to infinity as $i\to \infty$, where $s_2(i)=v_2\binom{2i}{i}$. Thus, $G_j\in \mathbb{Z}_2\langle x\rangle$. Moreover, differentiating gives
			\[
			G_j'(x)=-s_jr_j^3\left(1+\sum_{i\geq1}\binom{-\frac{3}{2}}{i}(2^hs_jt_jx)^i
			\right).
			\]
			Since $\binom{-3/2}{i}=(2i+1)\binom{-1/2}{i}$, we have, for every $i\ge1$,
			\[
			v_2\left(2^{hi}\binom{-3/2}{i}\right)
			=(h-2)i+s_2(i)\ge h-1.
			\]
			Thus, factoring $2^{h-1}$ out of the nonconstant terms, we obtain
			\begin{equation}\label{gprime}
				G_j'(x)=-s_jr_j^3+2^{h-1}E_j(x),\qquad E_j\in\mathbb Z_2\langle x\rangle.
			\end{equation}
			
			The square root $r_j(1+2^hs_jt_j^{-1}x)^{-1/2}=2^{h-1}G_j+r_j\equiv r_j\pmod 8$, so the $c_{\gamma}\left(r_j(1+2^hs_jt_j^{-1}x)^{-1/2}\right)$ in Lemma~\ref{classical} are constant. Expanding each Kloosterman sum into two terms expresses $\mathcal C(d)$ as a linear combination of the four sums
			\[
			R_{\epsilon_1,\epsilon_2}=\sideset{}{^{\star}}\sum_{x\bmod2^u}e_{2^u}(H(x)),
			\qquad H(x)=\epsilon_1G_1(x)-\epsilon_2G_2(Y(x))
			\]
			with coefficients of absolute value $2q$. Therefore
			\begin{equation}\label{fourbranches}
				|\mathcal C(d)|\le8q\max|R_{\epsilon_1,\epsilon_2}|.
			\end{equation}
			Since $d$ is even and $b_1,b_2$ are odd, we have $Y,Y'\in\mathbb Z_2\langle x\rangle$. Applying the chain rule and \eqref{gprime}, then we have
			\begin{equation}\label{hprime}
				H'(x)=A+\frac{B}{(b_1-dx)^2}+2^{h-1}E(x),
				\quad A=-\epsilon_1s_1r_1^3, B=\epsilon_2s_2r_2^3b_1b_2,
			\end{equation}
			where $A,B$ are odd units and $E\in\Z_2\langle x\rangle$.
			
			If $u=1$, the sum $R$ has only one term and the first assertion follows. We now assume $u>1$. Let $n=\lfloor u/2\rfloor\ge1$, $j=\lceil u/2\rceil$, and write $x=\beta+2^jz$. Since $2j\ge u$, Taylor expansion at $\beta$ gives
			\[
			R=2^n\mathop{\sstar}_{\substack{\beta\bmod2^j\\H'(\beta)\equiv0(\bmod 2^n)}}e_{2^u}(H(\beta)).
			\]
			Let $N_n$ denote the number of solutions modulo $2^n$, then $|R|\le2^jN_n$.
			
			If $h\ge n+1$, the congruence equation is $A(b_1-dx)^2+B\equiv0\pmod{2^n}$. When $v\ge n$, the trivial bound suffices. Otherwise, put $y=b_1-dx$. We first show that $y^2\equiv-\overline{A}B\pmod{2^n}$
			has at most four solutions modulo $2^n$. Suppose that odd $y_0$ is one such solution. Then any
			solution $y$ satisfies $2^n\mid(y+y_0)(y-y_0)$. Since $(y+y_0)-(y-y_0)=2y_0$, exactly we have
			\[
			2\parallel y+y_0\qquad\text{or}\qquad2\parallel y-y_0.
			\]
			Consequently, the other is divisible by $2^{n-1}$. Therefore, the four solutions are precisely $y\equiv \pm y_0\pmod{2^n}$, and $y\equiv \pm y_0+2^{n-1}\pmod{2^n}$. For each $y$, the congruence
			$b_1-dx\equiv y\pmod{2^n}$ has a unique solution $x\bmod 2^{n-v}$, and hence gives $2^v$ solutions modulo $2^n$. Thus $N_n\le4\cdot2^v$.
			
			If $v\le h-3$, for two distinct odd $x,y\in\Z$, we have
			\[
			H'(x)-H'(y)=\frac{Bd(x-y)(2b_1-d(x+y))}{(b_1-dx)^2(b_1-dy)^2}+2^{h-1}(E(x)-E(y)).
			\]
			It follows immediately that
			\[
			v_2(H'(x)-H'(y))=v+1+v_2(x-y).
			\]
			When $n>v+1$, any two solutions agree modulo $2^{n-v-1}$, and hence gives $2^{v+1}$ solutions modulo $2^n$. In the remaining case, a direct count also gives $N_n\le2^{v+1}$. This proves \eqref{corrbound}.
			
			If $n\ge h$ and $2^{h-2}\mid d$. 
			\[
			A+B{b_1}^{-2}\equiv 0\pmod {2^{h-2}}.
			\]
			Above equation implies \eqref{corrsupport}. If no sign branch satisfies this condition, every $R$ vanishes.
		\end{proof}

		For odd $A$, define the periodic function
		\[
		F_A(x)=\begin{cases}S(1,A\bar x;q),&x\text{ odd},\\0,&x\text{ even}.\end{cases}
		\]
		\begin{lemma}\label{Fsum}
			For $q=2^\gamma$, $\gamma\ge2$, and every integer $t$,
			\begin{equation}\label{auto}
				\sum_{x\bmod q}F_A(x)\overline{F_A(x+t)}
				=q\,c_q(t)=\begin{cases}q^2/2,&q\mid t,\\-q^2/2,&q/2\mid t,\ q\nmid t,\\0,&q/2\nmid t.\end{cases}
			\end{equation}
		\end{lemma}
		\begin{proof}
			Let $\widehat F_A(r)=\sum_xF_A(x)e_q(-rx)$. For odd $r$, expand $S$ and substitute $y=-rx$ to obtain
			$$
			\begin{aligned}
				\widehat F_A(r)=\sideset{}{^{\star}}\sum_{y,z\bmod q}e_q(y+z-Ar\overline{yz})=q\Kl(-Ar;q).
			\end{aligned}
			$$
			For even $r$, the change of variables $x\mapsto x+q/2$ changes the sign of the summand, so $	\widehat F_A(r)$ vanishes. By Lemma~\ref{klformula}, we have $|\widehat F_A(r)|^2=q^2$ for odd $r$. Using Fourier inversion, the left-hand side of \eqref{auto} is
            \[
			\begin{aligned}
				&	\frac{1}{q^2}\sum_{r,s\bmod q}\widehat F_A(r)\overline{\widehat F_A(s)}e_q(-st)\sum_{x\bmod q}e_q((r-s)x)\\
				&\hskip2in=\frac{1}{q}\sum_{r\bmod q}|\widehat F_A(r)|^2e_q(-rt)=q\sideset{}{^{\star}}\sum_{r\bmod q}e_q(-rt)=qc_q(t).
			\end{aligned}
            \]
		\end{proof}
		\subsection{Proof of Theorem \ref{main2}}
		We now prove Theorem~\ref{main2}. By Lemma~\ref{klformula}, the trivial bound suffices when $M\le q^{2/3}$. We may therefore assume that $q^{2/3}\leq M\leq q^{11/10}$. Throughout the paper we shall take $C=q^{100}$, let $\mathscr P=[C,2C]\cap\{\text{primes}\}$ and $\mathcal P=|\mathscr P|$. Choose a small parameter $0<\delta<10^{-3}$, to be specified, and set $Z=q^{\delta}$.
		
		Let $W\in C_c^\infty((0,\infty))$ such that $W=1$ on $\supp V$. Let $R=2^r\le q^{4/5}$, to be chosen later. Since $Rc$ is much larger than $M$, then for each $c\in\mathscr P$ we have the identity
		\[
		T=\sum_{m,n\ge1}\lambda_f(m)\Kl(an;q)e_M(Z(m-n))\mathbf1_{Rc\mid m-n}V\left(\frac{m}{M}\right)W\left(\frac{n}{M}\right).
		\]
		Averaging over $c$ and detecting the condition $Rc\mid m-n$ by additive characters, we may write $T$ as
		\[
		\frac{1}{R\mathcal P}\sum_{c\in \mathscr P}\frac{1}{c}\sum_{m,n\ge1}\lambda_f(m)\Kl(an;q)e_M(Z(m-n))\sum_{b\bmod Rc}e_{Rc}(b(m-n))V\left(\frac{m}{M}\right)W\left(\frac{n}{M}\right).
		\]
		Now let $(b,R)=P$, and write $bP$ for $b$. Hence the above formula equals
		\[
		\frac{1}{R\mathcal P}\sum_{c\in \mathscr P}\frac{1}{c}\sum_{m,n\ge1}\lambda_f(m)\Kl(an;q)e_M(Z(m-n))\sum_{P\mid R}\sum_{b\bmod Pc}e_{Pc}(b(m-n))V\left(\frac{m}{M}\right)W\left(\frac{n}{M}\right),
		\]
		where $P$ runs through powers of $2$. Now we split the sum into two cases $(b,Pc)=1$ and $(b,Pc)>1$. For $(b,Pc)>1$, we can bound trivially and get $O(M^2/C)$. Hence we have 
		\begin{equation}\label{decomp}
			T=\sum_{P\mid R}T_P+O(M^2/C),
		\end{equation}
		where
		$$
		\begin{aligned}\label{tpstart}
			T_P={\frac{1}{R\mathcal P}}\sum_{c\in\mathscr P}{\frac{1}{c}}\ \
			\sideset{}{^{\star}}\sum_{b\bmod Pc}\sum_{m,n\ge1}\lambda_f(m)\Kl(an;q)e_{Pc}(b(m-n))e_M(Z(m-n))V\left(\frac{m}{M}\right)W\left(\frac{n}{M}\right).
		\end{aligned}
		$$
		We first apply the Voronoi summation formula to the $m$-sum. By Lemma~\ref{Vor1}, we obtain
		\[
		\sum_{m\geq1}\lambda_f(m)e_{Pc}(bm)e_M(Zm)V\left(\frac{m}{M}\right)=\frac{M}{Pc}	\sum_{m\geq1}\lambda_f(m)e_{Pc}(-m\bar{b})I_1(m,c),
		\]
		where
		\begin{equation}\label{I1}
			I_1(m,c)=2\pi i^\kappa\int_{0}^{\infty}V(x)e(Zx)J_{\kappa-1}\left(\frac{4\pi\sqrt{mMx}}{Pc}\right)\mathrm{d} x.
		\end{equation}
		Making the change of variables $n=\alpha+jq$, the $n$-sum is transformed into the following sum:
		\[
		\sum_{\alpha\bmod q}\Kl(a\alpha;q)e_{Pc}(-b\alpha)\sum_{j\geq1}e_{Pc}(-b(\alpha+jq))e_M(-Z(\alpha+jq))W\left(\frac{\alpha+jq}{M}\right).
		\]
		We next apply the Poisson summation to the $j$-sum and obtain
		\begin{equation}\label{poi}
			\sum_{j}e_{Pc}(-b(\alpha+jq))e_M(-Z(\alpha+jq))W\left(\frac{\alpha+jq}{M}\right)=\frac{M}{q}\sum_{j}e_{q}(j\alpha)\int_{\mathbb{R}}W(x)e\left(-Zx-\frac{xM}{qc}\left(\frac{q}{P}b+cj\right)\right)\mathrm{d} x.
		\end{equation}
		Expanding $\Kl$ and applying orthogonality to the $\alpha$-sum, we have
		\begin{equation}\label{alpha}
			\sum_{\alpha\bmod q}\Kl(a\alpha,q)e_q(j\alpha)=\sideset{}{^{\star}}\sum_{\substack{x,y\bmod q\\a\overline{xy}+j\equiv 0(\bmod q)}}e_q(x+y)=\sideset{}{^{\star}}\sum_{x\bmod q}e_q(x-a\overline{jx})=S(1,-a\bar{j};q).
		\end{equation}
		Making the substitution $(q/P)b+cj=-k$ and combining the sums in \eqref{poi} and \eqref{alpha}, we obtain
		\[
		\frac{M}{q}\sum_{\substack{k\\c\mid{k+b(q/P)}}}F_{ca}\left(k+b\frac{q}{P}\right)I_2(k,c),
		\]
		where
		\begin{equation}\label{I2}
			I_2(k,c)=\int_0^\infty W(x)e\left(-Zx+\frac{kMx}{qc}\right)\,dx.
		\end{equation}
		Now by the Chinese Remainder Theorem, together with $b\equiv-k\overline{q/P}\pmod c$, we may write $b\equiv-k\overline{q/P}P\bar{P}+b'c\bar{c}\pmod{Pc}$. Then we get
		\[
		e_{Pc}(-\bar{b}m)=e_P(-\overline{b'c}m)e_c(\overline{P^2}q\bar{k}m).
		\]
		
		The corresponding finite sum is
		\begin{equation}\label{charfactor}
			e_c(\overline{P}^{\,2}q\bar k m)
			\sideset{}{^{\star}}\sum_{b\bmod P}F_{ca}(k+b(q/P))e_P(-\overline{bc}m).
		\end{equation}
		If $k$ is even, then $k+(q/P)b$ is even and this term vanishes. When $P=1$, the sum modulo $1$ is interpreted as a single term. Consequently,
		\begin{equation}\label{tpexact}
			T_P=\frac{M^2}{qRP\mathcal P}
			\sum_{c\in\mathscr P}{\frac{1}{c^2}}\sum_{m\ge1}\lambda_f(m)
			\sum_{\substack{k\in\Z\ \mathrm{odd}\\(k,c)=1}}
			e_c(\overline P^{\,2}q\bar k m)\,\mathcal A_c(m,k)\,I_1(m,c)I_2(k,c),
		\end{equation}
		where $\mathcal A_c$ denotes the sum modulo $P$ in \eqref{charfactor}. We now proceed to estimate $T_p$ in the remainder of this subsection.

		For convenience, we define $J(m,k,c)=I_1(m,c)I_2(k,c)$ and set
		\[
		H=\frac{Z^2P^2C^2}{M},\qquad D=\frac{ZqC}{M}.
		\]
		Let $h(x)=-Zx+\frac{kM x}{qc}$. Then we have
		\[
		h'(x)=-Z+\frac{kM}{qc}, \quad h^{(i)}(x)=0\ (i\geq 2).
		\]
		Hence, by the stationary phase method (see Lemma \ref{stationary} (i)), we have $I_2(k,c)\ll q^{-2026}$ unless $k\asymp D$. 
		
		By Lemma 3.2 of \cite{LS2021}, for $x\gg1$, we have
		\[
		J_{\kappa-1}(x)=x^{-1/2}\left(c_0 e\left(\frac{x}{2\pi}\right)+d_0 e\left(-\frac{x}{2\pi}\right)\right)+O(x^{-3/2}).
		\]
		Thus, the integral in $I_1(m,c)$ can be written as
		\[
		\int_{0}^{\infty}V(x)e\left(Zx\pm\frac{2\sqrt{mMx}}{Pc}\right)\mathrm{d} x.
		\]
		Similarly, Lemma~\ref{stationary} (i) gives $I_1(m,c)\ll q^{-2026}$ unless $m\asymp H$. For $x\ll1$, combining the bound $x^j J_{\kappa-1}^{(j)}\ll1$ with repeated integration by parts, we have
		\[
		I_1(m,c)=\frac{2\pi i^\kappa}{(-2\pi i Z)^{j}}\int_0^\infty e(Zx)\frac{d^j}{dx^j}\left(W(x)J_{\kappa-1}\!\left(\frac{4\pi\sqrt{mMx}}{Pc}\right)\right)\mathrm{d} x\ll Z^{-j}.
		\]
		Thus this contribution is negligible. We absorb $(C/c)^2$ into $J$ and apply a smooth partition of unity to the $m$-variable. By the Cauchy--Schwarz inequality, we obtain
		\begin{equation}\label{cauchy}
			|T_P|\ll\frac{M^{3/2}Z}{qR\mathcal P C}\Omega_P^{1/2},
		\end{equation}
		where
		\[
		\Omega_P=\sum_{m\in\Z}V_1\left(\frac{m}{H}\right)
		\left|\sum_{c\in\mathscr P}\sum_{\substack{k\ \mathrm{odd}\\(k,c)=1}}
		e_c(\overline{P^2}q\bar k m)\mathcal A_c(m,k)J(m,k,c)\right|^2.
		\]
		We are therefore reduced to estimating $\Omega_P$.
        \begin{lem}\label{deriv}
        We have
        \begin{equation}
			\left|H^iD^j \frac{\partial^i}{\partial m^i}\frac{\partial^j}{\partial k^j}J(m,k,c)\right|\ll Z^{i+j}.
		\end{equation}
        \end{lem}
        \begin{proof}
        Put
        \[
        z=\frac{4\pi\sqrt{mMx}}{Pc}\asymp Z.
        \]
        By repeatedly applying the chain rule, we obtain
        \[
        \frac{\partial^i}{\partial m^i}I_1(m,c)\ll \left(\frac{\partial z}{\partial m}\right)^iJ_{\kappa-1}^{(i)}(z)\ll(Z/H)^i,
        \]
        since $J_{\kappa-1}^{(i)}(z)\ll1$ (see \cite[p.\ 336]{Topacogullari2016}). A similar argument applies to $I_2(k,c)$. This completes the proof.
        \end{proof}

        Opening the absolute square, we break the $m$-sum into congruence classes modulo $L=Pc_1c_2$, and then apply the Poisson summation formula to the sum over $m$. It is therefore sufficient to consider the following sum
        $$
        \begin{aligned}
        \frac{1}{L}\sum_{c_1,c_2,k_1,k_2}\ \sum_{r\in\mathbb{Z}}&\int_{\mathbb{R}}V_1\left(\frac{y}{H}\right)J(y,k_1,c_1)\overline{J(y,k_2,c_2)}e_{L}(-ry)\mathrm{d} y\\\times&\sum_{\beta\bmod L}e_{c_1}(\overline{P^2}q\bar k_1 \beta)e_{c_2}(-\overline{P^2}q\bar k_2 \beta)A_{c_1}(\beta,k_1)\overline{A_{c_2}(\beta,k_2)}e_{L}(r\beta).
        \end{aligned}
        $$
        By integrating parts repeatedly and applying Lemma~\ref{deriv}, we may write the integral factor as
        \[
        H\mathcal I(k_1,k_2,r;c_1,c_2):=
        H\int_{\mathbb{R}}V_1(y)J(Hy,k_1,c_1)\overline{J(Hy,k_2,c_2)}e_{L}(-rHy)\mathrm{d} y\ll H\left(\frac{LZ}{rH}\right)^A,
        \]
        for any $A\geq 1$. Since $\frac{LZ}{rH}\asymp\frac{M}{rZP}$, this allows us to truncate the $r$-sum at
        \begin{equation}\label{F}
			|r|\le F:=M/P
		\end{equation}
        with a negligible error term. We next consider the sum over $\beta$. By the Chinese remainder theorem, we may write $\beta \equiv \beta_1 P + \beta_2 c_1 c_2 \pmod{L}$ uniquely. The sum over $\beta$ then equals
        $$
        \begin{aligned}
        &\mathop{\sstar\sstar}_{b_1,b_2\bmod P}F_{c_1 a}\left(k_1+b_1\frac{q}{P}\right)\overline{F_{c_2 a}\left(k_2+b_2\frac{q}{P}\right)}\\&\times\sum_{\beta_1 \bmod c_1c_2}e_{c_1c_2}(\beta_1\bar{P}q(c_2\overline{k_1}-c_1\overline{k_2}+\overline{q/P}r))\sum_{\beta_2\bmod P}e_P(\beta_2(c_1\overline{b_2}-c_2\overline{b_1}+r)).
        \end{aligned}
        $$
        Using the orthogonality of characters in the sum over $\beta_1$ and $\beta_2$, we obtain
        \begin{equation}\label{congruence}
        c_1\overline{k_2}-c_2\overline{k_1}\equiv \overline{q/P}r\pmod{c_1c_2}\quad\text{and}\quad c_2\bar b_1-c_1\bar b_2\equiv r\pmod P.
        \end{equation}
        Thus $\Omega_P$ equals
        \begin{equation}\label{Omega}
			H\sum_{r\ll F}\mathop{\sum\sum}_{c_1,c_2\in\mathscr{P}}\mathop{\sum\sum}_{\substack{k_1,k_2\ \mathrm{odd}\\(k_1,c_1)=(k_2,c_2)=1\\c_1\overline{k_2}-c_2\overline{k_1}\equiv \overline{q/P}r(\bmod{c_1c_2})}} \mathfrak C_2\mathcal I,
		\end{equation}
        where
        \begin{equation}\label{c2def}
			\mathfrak C_2=\mathop{\sstar\sstar}_{\substack{b_1,b_2\bmod P\\c_2\bar b_1-c_1\bar b_2\equiv r(\bmod P)}}
			F_{c_1 a}\left(k_1+b_1\frac{q}{P}\right)\overline{F_{c_2 a}\left(k_2+b_2\frac{q}{P}\right)}.
		\end{equation}
        Now we split $\Omega_P$ into $\Omega_P(r=0)$ and $\Omega_P(r\neq0)$. For this purpose, we use the following lemmas.
        \begin{lem}
        If $r=0$, then we have
        \begin{equation}\label{eq0}
        \Omega_{P}(r=0)\ll \frac{q^2P^3Z^3C^4}{M^2}.
        \end{equation}
        \end{lem}
        \begin{proof}
        When $r=0$, condition $c_1\overline{k_2}-c_2\overline{k_1}\equiv \overline{q/P}r\pmod{c_1c_2}$ forces $c_1=c_2=c$ and $k_1\equiv k_2\pmod{c}$. We may write $k_2=k_1+c\ell$, where $\ell\ll D/C\asymp Zq/M$. The congruence condition modulo $P$ gives $b_1=b_2$, so
        \[
        \mathfrak C_2(k)=\sideset{}{^{\star}}\sum_{b\bmod P}
			F_{c a}\left(k+b\frac{q}{P}\right)\overline{F_{c a}\left(k+c\ell+b\frac{q}{P}\right)}.
        \]
        We now consider the $k$-sum and introduce a smooth partition of unity in the $k$-variable. We temporarily include the terms with $c \mid k$, denoted by $\Omega_{P,0}$. We then split the $k$-sum into residue classes modulo $q$ and apply the Poisson summation formula, getting
        \[
        \sum_{k}V_1\left(\frac{k}{D}\right)\mathfrak C_2(k)\mathcal I(k,k+c\ell,0;c,c)=\frac{1}{q}\sum_{k\bmod q}\mathfrak C_2(k)\sum_{s\in\mathbb{Z}}e_{q}(ks)\int_{\mathbb{R}}V_1\left(\frac{y}{D}\right)\mathcal I(y,y+c\ell,0;c,c)e_{q}(-sy)\mathrm{d} y.
        \]
        
        For $s=0$, by the Lemma \ref{Fsum}, the $k$-sum is bounded by $O(\varphi(P)\,q\,c_q(c\ell))$. Since $c_q(c\ell)$ vanishes unless $q/2 \mid \ell$, only $O(1+Z/M)$ integers $\ell$ with $|\ell| \ll Zq/M$ contribute. It follows from \eqref{Omega} that
        \begin{equation}\label{Omega0}
        \Omega_{P,0}\ll C H \frac{D}{q}\varphi(P)q^2\asymp\frac{q^2P^3Z^3C^4}{M^2}.
        \end{equation}
        For $s\neq 0$, by integrating by parts repeatedly and applying Lemma~\ref{deriv}, we have
        $$
        \begin{aligned}
        &\int_{\mathbb{R}}V_1(y)\mathcal I(Dy,Dy+c\ell,0;c,c)e_{q}(-sDy)\mathrm{d} y\\&\hskip1in\ll\left(\frac{q}{sD}\right)^A\mathcal I^{(A)}(Dy,Dy+c\ell,0;c,c)(y)\\
        &\hskip1in\ll\left(\frac{q}{sD}\right)^A\frac{d^A}{dy^A}\int_{\mathbb{R}}V_1(x)J(Hx,Dy,c)\overline{J(Hx,Dy+cl,c)}e_{L}(-rHx)\mathrm{d} x\\
        &\hskip1in\ll\left(\frac{qZ}{sD}\right)^A
        \end{aligned}
        $$
        for any $A\geq 1$. Since $\frac{qZ}{sD}=\frac{M}{sC}$ and $C=q^{100}$ is sufficiently large, the resulting contribution is negligible.

        We now subtract the terms with $c \mid k$ that we temporarily included. There are at most $O(1+D/C)$ possibilities each for $k$ and $\ell$, and $\mathfrak C_2(k)$ can be bounded by $O(Pq)$ trivially. Their total contribution is
        \[
        \ll CH\left(1+\frac{D}{C}\right)^2\frac{D}{q}Pq,
        \]
        which is smaller than the right-hand side of \eqref{Omega0} since $M\geq q^{2/3}$.
        \end{proof}
        
        We now consider the case $r \neq 0$. 
        \begin{lem}
         We have
        \begin{equation}\label{neq0}
        \Omega_{P}(r\neq0)\ll Z^4P^{3/2}C^4q(1+q/M)^2.
        \end{equation}
        \end{lem}
        \begin{proof}
        Suppose that $P=2^u\ge2$. Applying the first case of Lemma \ref{local}, we obtain
        \[
		|\mathfrak C_2|\ll qP^{1/2}2^{v_2(r)}.
		\]
        The condition $c_1\overline{k_2}-c_2\overline{k_1}\equiv \overline{q/P}r\pmod{c_1c_2}$ implies $k_1\equiv-c_2(q/P)\bar{r}\pmod{c_1}$, since $c_1>|r|$. A similar argument applies to $k_2$. Hence, there are at most $O((1+D/C)^2)$ pairs $(k_1,k_2)$. Moreover,
        \[
		\sum_{0<|r|\ll F}2^{v_2(r)}\ll F\log(2F).
        \]
        Therefore the contribution from the cases $h\ge\lfloor u/2\rfloor+1$ or $v\le h-3$ is bounded by
		\begin{equation}\label{Bdef}
			HC^2P^{1/2}qF(1+D/C)^2\ll Z^4P^{3/2}C^4q(1+q/M)^2.
		\end{equation}

		In the remaining cases, $2^{h-2}\mid d$. Put
		\[
		Q=2^{h-2}=q/(4P),\quad d=Q\lambda,\quad
		0<|\lambda|\ll\Lambda:=4M/q.
		\]
        We also have the trivial bound $|\mathfrak C_2|\ll qP$. Applying the second case of Lemma \ref{local}, we have
        \[
        \epsilon_1\overline{c_1a}r_1^3c_2\equiv\epsilon_2\overline{c_2a}r_2^3c_1\pmod Q\,
        \]
        where $r_j^2=k_j^{-1}c_j a$. Taking square on both sides we have
        \begin{equation}\label{cubiccondition}
			c_1k_1^3\equiv c_2k_2^3\pmod Q.
		\end{equation}

		Reducing the first congruence in \eqref{congruence} modulo $c_1$ and $c_2$ in turn, we get
		\begin{equation}\label{linearcorrect}
			4c_2=\delta_1c_1+m\lambda,\qquad
			4c_1=\delta_2c_2-n\lambda.
		\end{equation}
		The two terms involving $\lambda$ must have opposite signs. Since $D\Lambda/C\ll Z$, we have $\delta_j\ll Z$.
		
		If $\delta_1\delta_2\ne16$, fixing $k_1,k_2,\lambda,\delta_1,\delta_2$ determines at most one pair $(c_1,c_2)$ through the linear system. The number of possible choices of $(c_1,c_2,k_1,k_2,\lambda,\delta_1,\delta_2)$ is therefore $\ll D^2\Lambda Z^2$.
		
		If $\delta_1\delta_2=16$, elimination gives $\delta_2k_1=4k_2$. Since $k_1,k_2$ are odd, $v_2(\delta_2)=2$, which yields
		\[
		\delta_1=\delta_2=4,\ k_2=k_1,
		\qquad\text{or}\qquad
		\delta_1=\delta_2=-4,\ k_2=-k_1.
		\]
		Substituting these two cases into \eqref{linearcorrect} and \eqref{cubiccondition}, respectively, gives
		$4(c_2\mp c_1)=k_1\lambda$ and $Q\mid(c_2\mp c_1)$. Consequently,
		\begin{equation}\label{sparse}
			4Q\mid\lambda.
		\end{equation}
		There are at most $O(\Lambda/Q)$ such nonzero $\lambda$. Fixing $k_1,\lambda,c_1$ determines the other variables. Thus, the number of $(c_1,c_2,k_1,k_2,\lambda)$ is $\ll DC\Lambda/Q$. Consequently, the contribution from the case $h\le\lfloor u/2\rfloor$ and $v\ge h-2$ can be bounded by
		\begin{equation}\label{degenerate}
			HqP\left(D^2\Lambda Z^2+DC\Lambda/Q\right).
		\end{equation}
        A simple calculation shows that if
		\begin{equation}\label{Pcondition}
			P\le\min\{q^{2/3}(1+M/q)^{4/3},q^{4/5}\},
		\end{equation}
		then the above contribution is also controlled by \eqref{Bdef}.
	
		Finally, we consider the case $P=1$, and we do not use Lemma~\ref{local} here. Expression \eqref{c2def} is simply a product of two Kloosterman sums and is $O(q)$. The outer congruence still gives $O((1+D/C)^2)$ pairs $(k_1,k_2)$. The length of $r$ has size $O(M)$, so
		\[
		\Omega_{1,\ne0}\ll HMC^2q(1+D/C)^2\ll Z^4qC^4(1+q/M)^2,
		\]
		which can be bounded by \eqref{Bdef}.
		\end{proof}
        
		Combining \eqref{eq0}, \eqref{neq0} and \eqref{cauchy}, we obtain
		\[
		|T_P|\ll \frac{M^{1/2}P^{3/2}Z^{5/2}}{R}
		+\frac{q^{1/2}M^{1/2}P^{3/4}Z^{3}}{R}(1+M/q).
		\]
		Choosing $\delta>0$ sufficiently small in terms of $\varepsilon$, any power of $Z$ can be absorbed into $q^{\varepsilon}$. Using \eqref{decomp}, we obtain
		\begin{equation}\label{beforeopt}
			|T|\ll q^{\varepsilon}M^{1/2}\{R^{1/2}+q^{1/2}R^{-1/4}(1+M/q)\}.
		\end{equation}
		Choose
		\[
		R=2^{\left\lfloor\log_2(q^{2/3}(1+M/q)^{4/3})\right\rfloor}.
		\]
		By \eqref{range}, $R\le q^{4/5}$. If $\gamma\ge20$, every $P=2^u\le R$ satisfies $\gamma-u\ge\gamma/5\ge4$. Condition \eqref{Pcondition} also holds. Substituting $R$ into\eqref{beforeopt} yields
		\[
		|T|\preceq q^{1/3}M^{1/2}(1+M/q)^{2/3}.
		\]
		This proves \eqref{strong}.
		
        If $2\le\gamma<20$, there are only finitely many possible $q$, and \eqref{range} also bounds $M$. Thus, the contribution from these cases can be absorbed into the implied constant. The theorem follows.

        \subsection*{Acknowledgements}
        The author thanks Prof. Jianya Liu for his encouragement. She is especially grateful to Prof. Yongxiao Lin for his helpful suggestions.


\begin{thebibliography}{9}
			
			
			
			\bibitem{BFKMM2017}
			V.~Blomer, \'E.~Fouvry, E.~Kowalski, P.~Michel and D.~Mili\'cevi\'c, On moments of twisted $L$-functions, {\it Amer. J. Math.} 139 (2017), no.~3, 707--768.

            \bibitem{BM2015}
            V.~Blomer and D.~Mili\'cevi\'c, The second moment of twisted modular $L$-functions, {\it Geom. Funct. Anal.} 25 (2015), no.~2, 453--516.

            \bibitem{CFKRS2005}
            J.~B. Conrey, D.~W. Farmer, J.~P. Keating, M.~O. Rubinstein and N.~C. Snaith, Integral moments of $L$-functions, {\it Proc. London Math. Soc.} (3) 91 (2005), no.~1, 33--104.
            
			\bibitem{DK2015}
			S.~Das and R.~Khan, Simultaneous nonvanishing of Dirichlet $L$-functions and twists of Hecke--Maass $L$-functions, {\it J. Ramanujan Math. Soc.} 30 (2015), no.~3, 237--250.
			
			\bibitem{HB1981}
			D.~R. Heath-Brown, The fourth power mean of Dirichlet's $L$-functions, {\it Analysis} 1 (1981), no.~1, 25--32.

            \bibitem{IK2004}
            H.~Iwaniec and E.~Kowalski, {\it Analytic number theory}, American Mathematical Society Colloquium Publications, 53, Amer. Math. Soc., Providence, RI, 2004.

            \bibitem{IS1999}
            H.~Iwaniec and P.~Sarnak, Dirichlet $L$-functions at the central point, in {\it Number theory in progress, Vol.~2 (Zakopane-Ko\'scielisko, 1997)}, 941--952, de Gruyter, Berlin, 1999.

            \bibitem{KZ2023}
            R.~Khan and Z.~Zhang, On the error term in a mixed moment of $L$-functions, {\it Mathematika} 69 (2023), no.~3, 573--583.

            \bibitem{KS2003}
            H.~Kim and P.~Sarnak, Appendix~2 in Functoriality for the exterior square of $\mathrm{GL}_4$ and the symmetric fourth of $\mathrm{GL}_2$, {\it J. Amer. Math. Soc.} 16 (2003), 139--183.

            \bibitem{KPY2019}
            E.~M.~K{\i}ral, I.~Petrow and M.~P.~Young, Oscillatory integrals with uniformity in parameters, {\it J. Th\'eor. Nombres Bordeaux} 31 (2019), 145--159.
  
            \bibitem{KMV2002}
            E.~Kowalski, P.~Michel and J.~VanderKam, Rankin--Selberg $L$-functions in the level aspect, {\it Duke Math. J.} 114 (2002), 123--191.
            
			\bibitem{LS2021}
			Y.~Lin and Q.~Sun, Analytic twists of ${\rm GL}_3\times{\rm GL}_2$ automorphic forms, {\it Int. Math. Res. Not. IMRN} (2021), no.~19, 15143--15208.
			
			\bibitem{Liu2015}
			S.-C.~Liu, Simultaneous nonvanishing of automorphic $L$-functions, {\it J. Number Theory} 147 (2015), 620--632.

            \bibitem{Milicevic2016}
            D.~Mili\'cevi\'c, Sub-Weyl subconvexity for Dirichlet $L$-functions to prime power moduli, {\it Compos. Math.} 152 (2016), no.~4, 825--875.

            \bibitem{MQW2025}
            D.~Mili\'cevi\'c, X.~Qin and X.~Wu, Bilinear forms with Kloosterman sums and moments of twisted $L$-functions, arXiv:2511.07550 (2025).

            \bibitem{Pascadi2026}
            A.~Pascadi, Non-abelian amplification and bilinear forms with Kloosterman sums, {\it Geom. Funct. Anal.} 36 (2026), no.~4, 1174--1243.

            \bibitem{Sharma2024}
            P.~K.~Sharma, Bilinear sums with $GL(2)$ coefficients and the exponent of distribution of $d_3$, {\it Proc. Lond. Math. Soc.} (3) 128 (2024), no.~3, Paper No. e12589, 42 pp.

            \bibitem{Shparlinski2019}
            I.~Shparlinski, On sums of Kloosterman and Gauss sums, {\it Trans. Amer. Math. Soc.} 371 (2019), no.~12, 8679--8697.

            \bibitem{Sono2019}
            K.~Sono, Simultaneous nonvanishing of Dirichlet $L$-functions and twists of Hecke--Maass $L$-functions in the critical strip, {\it Ann. Acad. Sci. Fenn. Math.} 44 (2019), no.~2, 1131--1157.
             
			\bibitem{Sound2007}
			K.~Soundararajan, The fourth moment of Dirichlet $L$-functions, in {\it Analytic number theory}, 239--246, Clay Math. Proc. 7, Amer. Math. Soc., Providence, RI, 2007.

            \bibitem{TW2025}
            Z.~Tang and X.~Wu, Mixed moments of twisted $L$-functions, arXiv:2512.09203 (2025).

            \bibitem{Topacogullari2016}
            B.~Topacogullari, The shifted convolution of divisor functions, {\it Q. J. Math.} 67 (2016), no.~2, 331--363.

            \bibitem{WZ2026}
            H.~Wang and T.~Zhu, Triple divisor type functions in arithmetic progressions to prime power moduli, arXiv:2609.16700 (2026).
            
			\bibitem{Wu2023}
			X.~Wu, The fourth moment of Dirichlet $L$-functions at the central value, {\it Math. Ann.} 387 (2023), no.~3--4, 1199--1248.
			
			\bibitem{Young2011}
			M.~P. Young, The fourth moment of Dirichlet $L$-functions, {\it Ann. of Math.} (2) 173 (2011), no.~1, 1--50.	
			
			
		\end{thebibliography}
	\end{document}